\documentclass[12pt]{amsart}
\usepackage{amssymb,latexsym,epsfig,color, hyperref}

\newcommand{\Z}{{\mathbb Z}}

\newcommand{\red}{\color{red}}

\usepackage{extarrows, url}

\newtheorem{theorem}{Theorem} 

\theoremstyle{definition}

\title[Divisibility Tests Unified
]
{Divisibility Tests Unified: \\ Stacking the Trimmings for Sums
}
\author{Edwin O'Shea}
\thanks{After this article was submitted to {\em Math Magazine}, 
Eric L. McDowell independently arrived at 
many of the same results in 
``\href{https://www.maa.org/press/periodicals/convergence/divisibility-tests-a-history-and-users-guide}
{Divisibility Tests: A History and User's Guide}'' 
which appeared in May 2018 in \textit{MAA Convergence}. What we call ``stacking'' is 
what that paper refers to as ``flowing''. McDowell's presentation deserves a 
wide readership and contains many references not addressed in this present paper.}
\address{Department of Mathematics \& Statistics, James Madison University, Harrisonburg, VA 22807-1911}
\email{osheaem@jmu.edu}
\date{\today}

\begin{document}

\begin{abstract}
Divisibility tests are algorithms that can quickly decide if one integer is divisible by another. 
There are many tests but most are either of the {\em trimming} or {\em summing} 
variety. Our goals are to present Zbikowski's family of trimming tests as one test and to unify the 
trimming and summing tests. We do the latter by showing, first, that the most effective summing tests, 
due to Khare, can be derived directly from the Zbikowski's test and, second, that the best known 
summing tests - the binomial tests - can be derived from an adapted form of Zbikowski's tests.  
We introduce the notion of {\em stacking}, the claim that a six year old 
would always choose 10 pennies over a dime, and use only basic divisibility properties to 
achieve our goals.

\end{abstract}

\maketitle

\subsection*{Introduction} 
The most well-known divisibility tests are the last digits tests for $2$ and $5$, 
the sum of digits test for $9$, and the alternating digit sum for $11$, 
but the oldest divisibility test is one for deciding divisibility by $7$. That test is 
at least fifteen hundred years old and is prescribed in the Talmud \cite[Abodah Zarah 9b]{Talmud} as follows: 
``If one does not know what the year is in the Sabbatical cycle of seven years, 
	let him... put aside the hundreds... and convert the remainder into Sabbatical Cycles 
	[of seven years each] after adding thereto two years for every complete century; 
	what is left over will give him the number of the given year in the current Sabbatical Cycle.''
In algebraic notation, the remainder when $7$ is divided 
into a given integer, written as $x+ 100 \cdot y$, equals that when $7$ is divided into $x + 2 \cdot y$. 
For example, to remainder when $7$ divides $32184 = 84 + 100 \cdot 321$ equals that when $7$ 
divides $84 + 2 \cdot 321 = 726$. 

The Talmud's test is the first of seventy or so listed in Dickson's encyclopedic 
{\em History of the Theory of Numbers} \cite[Chapter XII]{Dickson} and includes 
tests by luminaries such as Fibonacci, Lagrange, Pascal, and Sylvester. 
Tests that reinterpret those recorded by Dickson can be found in a number of relatively 
recent papers \cite{CheMou, Ganzell, Renault, Zazkis} and the sources referenced therein.
Among the tests is one for $7$ by Zbikowski \cite{Zbikowski} asserting that an integer $a$, written in 
the form $a = 10\bar{a} + a_0$, is divisible by $7$ if and only if $7$ divides the integer $\bar{a} - 2 a_0$.
For example, the test reduces $32184$ to $3218 - 2 \cdot 4 = 3210$, which can be applied 
again, reducing $3210$ to $321$, and again, reducing $321$ to $30$; since $7$ does not divide $30$ 
it does not divide $32184$. 

This {\em trimming} procedure, the given integer $a$ being ``trimmed'' to another with one digit less, 
is universally presented as being cut from a different cloth from the sum of digits tests for 
$9$ and $11$. We claim that this is not so by showing 
by deriving a family of summing tests, due to Khare \cite{Khare}, from Zbikowski's family of 
trimming tests. We can also show that the best known summing tests, the binomial tests, 
can also be derived from an adapted form of Zbikowski's tests. 
To the best of our knowledge this marriage of trimming and summing tests is new. 

In homage to the school venue where many of us were first exposed to divisibility tests, 
we will only require basic properties of the integers with a dash of the induction axiom; 
we will not use the binomial theorem or modular arithmetic. Our central tool is {\em stacking}, 
a decimal representation that is flexible enough to respect a six year old's choosing of 
ten pennies over one dime. 
The well known sum and alternating sum of digits tests for $9$ and $11$ follow as 
corollaries. We close with a brief comparative analysis of Khare's tests, the 
binomial summing tests and Zbikwoski's trimming tests, and how these tests in base 
$10$ generalize to any base. 

\subsection*{Divisibility Tests}

Rather than operate under a Justice Potter-like assumption \cite{Supreme}, that 
we all know a divisibility test when we see it, let us propose
a decent definition.
In most basic terms, a {\em divisibility test for an integer $q$} 
should be a function $f_q: \Z \rightarrow \Z$ 
such that $q$ divides $a$ if and only if $q$ divides $f_q(a)$ for every 
integer $a$. 
The identity function $f_q(a) = a$ is easy to compute but $q$ dividing $f_q(a)$ 
is no easier to decide than if $q$ divides $a$.  
The computation of the remainder in the classical division theorem, 
$f_q(a) = r$, might have the property that $q | r$ is easier to decide 
than $q | a$ but the computation of $f_q(a)$ is likely to be mentally 
difficult. We'd like to propose that a divisibility test $f_q(a)$ should 
be {\em easy} to compute and it ought to be {\em easier} to decide if 
$q$ divides $f_q(a)$ than if $q$ divides $a$. 
The terms ``easy'' and ``easier'' are ambiguous but one 
criterion for ``easy'' is that $f_q(a)$ is computable with relative ease.
``Easier'' could also mean a number of things but a desirable property might be 
that the number of digits in $f_q(a)$ is less than that in $a$. 
Note that any test $f_q$ is iterative, 
with $f^2_q(a) = f_q(f_q(a))$ being a test too for $q$ dividing $a$, and 
$f^3_q(a) = f_q(f_q(f_q(a)))$ too, etc..  

We promised to only use basic divisibility properties to derive our tests; no modular arithmetic or 
binomial theorem. To that effect, the following appear in number theory texts like Andrews \cite{And}. 

\vspace{.1cm}

\begin{tabular}{lll}
{\bf (1)} & \hspace{1in} & { If two integers $a$ and $s$ are both divisible by $q$ then} \\  
	  & \hspace{1in} & { their sum and difference, $a \pm s$ are also divisible by $q$.}  
\end{tabular}


\vspace{.1cm}

\begin{tabular}{lll}
{\bf (2)} & \hspace{1in} & { If $q$ is relatively prime to $10$ and $q$ divides $10 \cdot m$} \\ 
	    & \hspace{1in} & {  then $q$ divides $m$.}  
\end{tabular}

We can write an integer $a$ as 
$ 
a = a_n a_{n-1} \ldots a_2 a_1 a_0 = \sum_{k=0}^n 10^k \cdot  a_k,
$  
where each $0 \leq a_k \leq 9$. 
For shorthand, we denote the number of digits of $a$, $n+1$, as length$(a)$. 
Letting 
$
a_{[k,l]} :=  a_k a_{k-1}  \ldots  a_{l+1} a_l, 
$  
we can always write
$ 
a = 10^k \cdot a_{[n,k]} + a_{[k-1,0]}.  
$ 
As a special case, let $\bar{a}$ denote $a_{[n,1]}$ and write 
$$
a = 10 \cdot \bar{a} + a_0 \, \, \, \, \, \textup{and similarly,} \, \, \, \,  \, q = 10 \cdot \bar{q} + q_0.
$$
For example, if $a=32184$ then $a_4 = 3, a_3 = 2, a_2 = 1, a_1 = 8$ and $a_0 = 4$. The length of $a$ is $5$. We can write 
$32184$ in a variety of ways including $10^2 \cdot a_{[4,2]} + a_{[1,0]} = 10^2 \cdot 321 + 84$ and 
$10 \cdot \bar{a} + a_0 = 10 \cdot 3218 + 4$.

It is left as an exercise to apply claim {\bf (1)} to derive the last digit tests $f_{2}(a) = f_{5}(a) = a_0$. 
More generally, $f_{2^k}(a) = f_{5^k}(a) = a_{[k-1,0]}$ 
are divisibility tests for $q=2^k$ and $q=5^k$. For example, $8$ divides $a=32184$ because $8$ divides 
$a$'s last three digits, $f_{2^3}(32184) = 184$. 
With the above notation, the Talmud test is $\textup{Tal}_7(a) = 2 \cdot a_{[n,2]} + a_1a_0$. It too can be proved 
using claim {\bf (1)}: letting $s = 98 \cdot a_{[n,2]}$, $7$ divides $a = 100 \cdot a_{[n,2]} + a_1a_0$ 
if and only if it divides $a - s =  2 \cdot a_{[n,2]} + a_1a_0$.

\subsection*{Zbikowski's Trimming Tests as One Test}

Zbikowski's test for $7$ is $T_7(a) = \bar{a} - 2 \cdot a_0$.  
On an example like $T_7(32184) = 3210$ we see that $T_7$ 
takes a given $a$ and ``trims'' it to another integer of length 
one less than the original $a$. This motivates the following definition: 

\vspace{.1cm}

\noindent 
\begin{tabular}{lll}
{\bf (}{\bf Trimming}{\bf )} 
& \hspace{.2in} & 
A divisibility test $f_q$ is called a {\em trimming} test if the length \\  
& \hspace{.2in} & of $f_q(a)$ is one less than the length of $a$, for almost every $a$. 
\end{tabular}

\vspace{.1cm}

We say ``almost'' because if $a$ is already a single digit there is nothing to be done 
and there are instances, like $T_7(49) = -14$, 
where the test maps a two-digit number to another two-digit number. We leave it as 
an exercise to show that if length$(a) \geq 3$ then $T_7(a)$ has shorter length than $a$. 

Here's why $T_7$ works on our running example of $a=32184$.  
By claim {\bf (1)}, 
we can subtract any multiple of $7$ from $2184$ and the result will be 
divisible by $7$ if $32184$ itself is divisible by $7$, so choose a multiple of $7$ 
that when subtracted from $32184$ leaves a zero in the last digit. Clearly, $21$ 
times the last digit of $32184$, namely $21 \cdot 4$ will serve this role. 
The difference is $32184 - 21 \cdot 4 =  (10 \cdot 3218 + 4) - ( (20 + 1) \cdot 4)$.  
The $4$'s cancel leaving a multiple of $10$. 
By claim {\bf (2)}, we can trim that right-most zero from 
$32184 - 21 \cdot 4 = 32100$ to get $3210$ and our decision of whether $7$ divides $32184$ 
becomes equivalent to deciding if $7$ divides $3210$. 

Zbikowski \cite{Zbikowski} extended this argument for every $a$ and for any $q$ with last 
digit equal to $1,3,7,$ or $9$. These tests have received considerable attention in 
recent papers by Zazkis \cite{Zazkis}, Cherniavsky and Mouftakhov~\cite{CheMou}, 
and Ganzell~\cite{Ganzell} and the reader can see a derivation of these tests 
there. We will not derive these tests here but wish to recast these tests as one test. 
First, Zbikowski's tests as four different cases, followed by examples. 
\begin{theorem} \label{the:trim} \textup{(Zbikowski \cite{Zbikowski})}
For every $q$ with last digit equal to either $1,3,7,$ or $9$, there is a trimming test $T_q(a)$ given by 
the following table.  

\begin{center}
\begin{tabular}{|c|c|c|c|c|}                                                        
\hline                                                                          
$q_0$ 	 & $1$ 			   & $3$ 			  & $7$ 			 & $9$ \\ \hline 
$c_q$    & $1$ 			   & $-3$ 			  & $3$ 			 & $-1$ \\ \hline 
$T_q(a) = t_{q}(a, c_q)$ & 
$\bar{a} - \bar{q} a_0$ & 
$\bar{a} + (3 \bar{q}+1)a_0$ & 
$\bar{a} - (3\bar{q}+2) a_0$ & 
$\bar{a} + (\bar{q}+1)a_0$ \\ \hline
\end{tabular}                                                                   
\end{center}

\end{theorem}

\begin{itemize}
\item[] If $q = 21$  then $\bar{q} = 2$ and $T_{21}(a)  \, = \, \bar{a} -  2 \cdot a_0 \, = \, \bar{a} - 2 a_0$.
\item[] \hspace{.1in} For $a=32184$,  $T_{21}(32184) = 3218 - 2 \cdot 4 = 3210$ and $T_{21}(3210) = 321.$ 
\item[] If $q = 13$  then $\bar{q} = 1$  and  $T_{13}(a)  \, = \, \bar{a} + (3 \cdot 1 + 1)a_0 \, = \, \bar{a} + 4 a_0$ . 
\item[] \hspace{.1in} For $a=32184$,   $T_{13}(32184) = 3218 + 4 \cdot 4 = 3234$ and $T_{13}(3234) = 339.$ 
\item[] If $q = 17$  then $\bar{q} = 1$  and  $T_{17}(a)  \, = \, \bar{a} - (3 \cdot 1 + 2)a_0 \, = \, \bar{a} -5 a_0$ . 
\item[] \hspace{.1in} For $a=32184$, $T_{17}(2184) = 3218 - 5 \cdot 4 = 3198$ and $T_{17}(3198) = 279.$ 
\item[] If $q = 39$  then $\bar{q} = 3$  and  $T_{39}(a)  \, = \, \bar{a} + ( 3+ 1 )a_0 \, = \, \bar{a} + 4 a_0$. 
\item[] \hspace{.1in} For $a=32184$, $T_{39}(32184) = 3218 + 4 \cdot 4 = 3234$ and $T_{39}(3234) = 339.$ 
\end{itemize}                                                                   

As expected the above examples trim one integer per iteration. The examples 
are for $q$'s with two digits but $q$ can be of any length, like $T_{181} = \bar{a} - 18 \cdot a_0$. 

Absent from previous expositions on Zbikowski is that the four tests reduce to one. 
First, it appears that  $T_{13} = T_{39}$ and $T_7 = T_{21}$. Using the table above, one can 
show that: 

\vspace{.1cm}

\begin{tabular}{lll}
{\bf (3)} & \hspace{.1in} & 
If  $q_0=3$  or  $7$  then  $\displaystyle T_q(a) = T_{3q}(a)$.  
\end{tabular}

\vspace{.1cm}

\noindent This reduces our four tests to only two, those $T_q$'s for which $q_0 =1$ or $q_0 = 9$.  
With $\left[  x \right]$ denoting the nearest integer to $x$ we leave it to the reader, using the 
table above, to confirm:

\vspace{.1cm}

\begin{tabular}{lll}
{\bf (4)} & \hspace{-.05in} & 
$
T_q(a) \, = \, \bar{a} + \omega_q \cdot a_0 \, \, \, \, \textup{where} \, \, \, \, 
\omega_{q} \, = \, 
\left\{ 
\begin{array}{rcl}
-\left[  \frac{q}{10} \right]  & \textup{if} & q_0 = 1 \\ 
\left[  \frac{q}{10} \right]  & \textup{if} & q_0 = 9
\end{array}
\right\} .
$
\end{tabular}

\vspace{.1cm}

In summary, Zbikowksi's test reads easily as one test: 
{\em If an odd divisor $q$ ends in $1$ or $9$ then divide $q$ by 
$10$ and round the result to the nearest integer; attach a sign of minus or plus 
to the result depending on whether you have rounded down or up for the signed weight $\omega_q$.  
If $q$ ends in $3$ or $7$ then triple $q$ and do as before; that is, $\omega_q = \omega_{3q}$. Zbikowski's 
test for $q$ dividing $a$ is then everything but the last digit of $a$ plus the signed weight $\omega_q$ 
times the last digit of $a$.} 

For example, to write a divisibility test for $q=17$ we triple $17$ to get $51$. For the signed weight 
$\omega_{17}$, divide $51$ by $10$ and round to the nearest integer to produce $5$; since we rounded 
down the signed weight must be negative and so $-5$ is the weight for the test for $q=17$. 
That is, $T_{17}(a) = \bar{a} - 5 a_0$. Likewise, $T_{79} = \bar{a} + 8 a_0$ since $79/10$ rounds to $8$ 
and the weight is positive since we rounded up (not down) to $8$. 

Using Zbikowski's trimming test $T_q(a) \, = \, \bar{a} + \omega_q \cdot a_0$ we shall derive Khare's 
general weighted sum of digits tests \cite{Khare}. Khare's summing tests $S_q$ match the usual tests 
for $9$ and $11$ but differ from the better known binomial tests for all other $q$. Nonetheless, we can 
also derive the usual binomial tests by adapting Zbikowski's tests to trim from the left rather than the 
right. This is all achieved by a form of child's play we call stacking.  

\subsection*{Stacking: Preferring Pennies to Dimes}

The trimming tests $T_9(a) = \bar{a} + a_0$ and $T_{11}(a) = \bar{a} - a_0$ are not the same yet 
look similar to the sum and alternating sum of digits tests respectively. These sum of digits tests are usually 
verified by modular arithmetic -- geometric series suffice too --  
but the trimming tests have only used the basic divisibility properties {\bf (1)} and {\bf (2)}. 
From the trimming tests $T_q$ we will derive the usual tests for $9$ and $11$ 
and Khare's {\em summing} tests for every $q$. We should first define what we mean by a summing test. 

\vspace{.1cm}

\noindent 
\begin{tabular}{lll}
{\bf (}{\bf Summing}{\bf )} & \hspace{.2in} & A divisibility test $f_q(a) = \sum_{j=0}^n \gamma_j a_{n-j}$   \\
	        			       & \hspace{.2in} & 
				       is called a {\em summing} test for $q$ if each $\gamma_j \in \Z$.  
\end{tabular}

\vspace{.1cm}

Let's investigate the trimming test $T_9(a) = \bar{a} + a_0$ with our running example $a=32184$ and see 
if we can get some ideas on how to derive the sum of digits test $3+2+1+8+4$. 
The trimming test applied iteratively is 

$$ 
32184 
\, {\overset{T_9}\longrightarrow} \, 3218 + 4 \, = \, 3222 
\, { \overset{T_9}\longrightarrow} \, 322+2 \, = \, 324 
\, {\overset{T_9}\longrightarrow} \, 32+4 \, = \, 36  
\, {\overset{T_9}\longrightarrow} \, 3+6 .  
$$

The summing and recursive trimming tests yield a different final output. 
We claim that they are equal provided that a ``stacking'' 
procedure intervenes. To explain the main idea, let's start with a 
non-trivial theorem, that of 
{\em every positive integer has a unique base $10$ representation}. This is mathematically 
respected but colloquially malleable. When writing checks we are allowed to express $1562$ in unambiguous 
but different ways, as both``one thousand, five hundred and sixty two'' and as 
``fifteen hundred and sixty two.''  The former is in keeping with strict mathematical 
practice yet the latter is customary even though $15$, the coefficient 
(allowing ourselves to call it that) of one hundred in the latter is not between $0$ and $9$. 

In the same vein, when adults add two integers, like $3218+4 = 3222$ 
that result from $T_9(32184)$, we simplify in concordance with unique representability. 
Computing the sum $3218+4$ is equivalent to giving an adult $3218$ cents as 
$321$ dimes and $8$ pennies and giving them a further $4$ pennies, with which the adult opts 
to exchange $8+4 = 12$ pennies for $1$ dime and $2$ pennies for a total of $322$ dimes and $2$ pennies. 
We are raised to value efficiency; the fewer coins, the better. 
However, given the same choice, a six-year old may opt to keep the $12$ pennies. 
She knows that $10$ pennies and $1$ dime both equal $10$ cents but $10$ pennies are far 
more fun to play with and easier to share than a dime and so she chooses to stack the pennies together. 
In other words, she might opt for $321$ dimes and $8+4 = 12$ pennies, that is 
$3218+4 = 10 \cdot 321 + 8 + 4 = 10 \cdot 321 + (8 + 4)$. 
Depending on her mathematical formalism, she would define stacking the pennies as follows.

\vspace{.1cm}

\noindent 
\begin{tabular}{lll}
{\bf (}{\bf Stacking}{\bf )} & \hspace{.2in} & Given an integer $r = 10\bar{r}+r_0$ and a (possibly empty) sum \\  
			     & \hspace{.2in} & of single-digit integers $s$ write the {\em stacking} of their sum \\  
			     & \hspace{.2in} &   
			     $r + {\red s} \, \, {\red {\xlongequal{\textup{Stack}}} } \, \,  10 \bar{r} + {\red (r_0 + s)}$. 			     
\end{tabular}

\vspace{.1cm}

For short, we write the stacking of $r$ and $s$ as Stack$(r+s)$. 
For example, stacking $3218$ and $4$ together 
equals the representation Stack$(3218+4) = 10 \cdot 321 + (8+4)$. 
Since stacking is nothing more than an alternative representation of $r+s$, $q$ divides 
$r+s$ if and only if $q$ divides Stack$(r+s)$. 

\subsection*{Stacking Zbikowski Trimmings for Khare's Summing Tests}

With stacking in mind, let's iteratively trim as before with $T_9$ but now follow each trimming 
with a stacking. 

\vspace{.1in}

$3218{\red 4} \, { \xlongrightarrow{\, \, \, \, T_9 \, \, \, \, }} \, \, 3218 + {\red 4} \, 
{\red {\xlongequal{\bf Stack}} } \, 10 \cdot 321 +{\red (8+4)}$

\hspace{.35in}
$
{ \xlongrightarrow{\, \, \, \, T_9 \, \, \, \, }} \, \, 
\, 321 + {\red (8+4)} \, 
{\red {\xlongequal{\bf Stack}} } \, 10 \cdot 32 +{\red (1+8+4)}$

\hspace{.35in}
$
{ \xlongrightarrow{\, \, \, \, T_9 \, \, \, \, }} \, \, 
32 + {\red (1 +8 + 4)} \, 
{\red {\xlongequal{\bf Stack}} } \, 10 \cdot 3 + {\red (2+1+8+4)}$

\hspace{.35in}
$
{ \xlongrightarrow{\, \, \, \, T_9 \, \, \, \, }} \, \, 
3 + {\red (2+ 1 +8 + 4)} \, 
{\red {\xlongequal{\bf Stack}} } \, {\red (3+2+1+8+4).}$

\vspace{.1in}

The above says that $(\textup{Stack} \circ T_9)^{4}(32184) = {\red (3+2+1+8+4)} = 18 =: S_9(32184)$, 
where the latter denotes the usual sum of the digits test for $9$.  
Let us see if iteratively trimming and stacking with $T_7(a) = \bar{a} + (-2)a_0$ can provide a sum-like 
test for $q=7$ using our running example $a=32184$.  

\vspace{.1in}

$
3218{\red 4} \, { \xlongrightarrow{\, \, \, \, T_7 \, \, \, \, }} \, \, 3218 + {\red (-2) \cdot 4} \, 
$

\hspace{.35in}
$
{\red {\xlongequal{\bf Stack}} } \, \,  10 \cdot 321 +{\red (8+ (-2) \cdot 4) }
$

\hspace{.35in}
$
{ \xlongrightarrow{\, \, \, \, T_7 \, \, \, \, }} \, \, 
\,  321 +{\red  (-2) \cdot (8+ (-2) \cdot 4)  } \, 
$

\hspace{.35in}
$
{\red {\xlongequal{\bf Stack}} } \, \, 
10 \cdot 32 +{\red (1 +  (-2) \cdot (8+ (-2) \cdot 4) ) }
$

\hspace{.35in}
$
{ \xlongrightarrow{\, \, \, \, T_7 \, \, \, \, }} \, \, 
32 +{\red (-2) \cdot (1 +  (-2) \cdot (8+ (-2) \cdot 4) ) } \, 
$

\hspace{.35in}
$
{\red {\xlongequal{\bf Stack}} } \, \, 
10 \cdot 3 + {\red 2 + (-2) \cdot (1 +  (-2) \cdot (8+ (-2) \cdot 4) ) }
$

\hspace{.35in}
$
{ \xlongrightarrow{\, \, \, \, T_7 \, \, \, \, }} \, \, 
3 + {\red (-2)(2 + (-2) \cdot (1 +  (-2) \cdot (8+ (-2) \cdot 4) ) )}
$

\hspace{.35in}
$
{\red {\xlongequal{\bf Stack}} } \, \, 
 {\red 3 + (-2)(2 + (-2) \cdot (1 +  (-2) \cdot (8+ (-2) \cdot 4) ) ).}
$

\vspace{.1in}

In other words,  
$(\textup{Stack} \circ T_7)^{4}(32184) \, = \, 
3 + (-2)^1 2 + (-2)^2 \cdot 1 + (-2)^3 \cdot 8 + (-2)^4 \cdot 4$. 
The above examples for $q=7$ and $q=9$ with $a=32184$ suggest 
summing tests with $\gamma_j = (-2)^j = \omega_7^j$ and $\gamma_j = 1 = \omega_9^j$ 
for $7$ and $9$ respectively. We claim this holds in general. 

\begin{theorem} \label{the:sum}  
If $T_q = \bar{a} + \omega_q a_0$ is a trimming test for $q$ then  
$
S_q(a) 
\, := \,  
\sum_{j=0}^{n} \omega_q^j a_{n-j}
$ 
is a summing test for $q$. 
\end{theorem}

The tests $S_q$ were presented in 1997 by Khare \cite{Khare} 
but their modular arithmetic proof does not involve trimming tests. 
Briefly, Khare's construction begins by choosing of $\gamma_q$ as the minimum residue representative 
of the inverse of $10$ modulo $q$. That is, $\gamma_q \equiv 10^{-1} \mod q$ of smallest size. 
Khare then proposes $S_q = \sum_{j=0}^{n} \gamma_q^j a_{n-j}$ is a test by virtue of  
$$ 
\gamma_q^n a \, \, = \, \, \sum_{j=0}^n \gamma_q^n {10^j} a_j \equiv S_q(a) \mod q. 
$$ 

It is straightforward to check that Khare's $\gamma_q$ equals Zbikowski's $\omega_q$. 
Our derivation of Khare's tests from Zbikowski's tests 
uses neither modular arithmetic or the binomial theorem and it unifies 
the trimming and summing families. 
Before proving the result, let's appreciate Khare's tests for some examples on $a=32184$: 

\begin{center}
\begin{itemize}
\item[] $S_{7}(32184) \, = \, 3+ (-2) \cdot 2 + (-2)^2 \cdot 1 + (-2)^3 \cdot 8 + (-2)^4 \cdot 4 = 3$
\item[] $S_{9}(32184) \, = \, 3 + 1 \cdot 2 + 1^2 \cdot 1 + 1^3 \cdot 8 + 1^4 \cdot 4 = 18$.
\item[] $S_{11}(32184) \, = \, 3 + (-1) \cdot 2 + (-1)^2 \cdot 1 + (-1)^3 \cdot 8 + (-1)^4 \cdot 4 = -2$.
\item[] $S_{17}(32184) \, = \, 3 + (-5) \cdot 2 + (-5)^2 \cdot 1 + (-5)^3 \cdot 8 + (-5)^4 \cdot 4 = 1518$.
\item[] $S_{39}(32184) \, = \, 3+ 4 \cdot 2 + 4^2 \cdot 1 + 4^3 \cdot 8 + 4^4 \cdot 4 = 1563$. 
\end{itemize}                                                                   
\end{center}

\bigskip 

\noindent {\em Proof of Theorem~\ref{the:sum} by Trimming and Stacking}. \,  
We will show, by induction on the length of $a$, that $S_q(a) = (\textup{Stack} \circ T_q)^{n}(a)$ 
whenever $a$ has length $n+1$. 

If $n=1$ then $a = a_1a_0$ has length two and 
$\textup{Stack}(T_q(a_1a_0)) = \textup{Stack}(\bar{a} + \omega_q a_0) = a_1 + \omega_q a_0$ 
as claimed. 
Assume that $S_q(a') = (\textup{Stack} \circ T_q)^{n-1}(a')$ 
for every $a'$ with length $n$ and consider any integer 
$a = a_n a_{n-1} \ldots a_2 a_1 a_0$ with length $n+1$. 
Applying $\textup{Stack} \circ T_q$ to this $a$ results in 
$
\textup{Stack}(T_q(a)) = \textup{Stack}(\bar{a} + \omega_q a_0)  
\, = \, 
10 a_{[n,2]} +(a_1 + \omega_q a_0), 
$ 
an integer with $n$ digits with last digit equal to $(a_1 + \omega_q a_0)$ 
to which the induction hypothesis applies; hence, 
\begin{align*}
(\textup{Stack} \circ T_q)^{n}(a) 
& =  \, (\textup{Stack} \circ T_q)^{n-1} ( \textup{Stack} \circ T_q(a))  \\
& = \, 
(\textup{Stack} \circ T_q)^{n-1} ( 10 a_{[n,2]} +(a_1 + \omega_q a_0)) \\ 
 & =  \, S_q( 10 a_{[n,2]} +(a_1 + \omega_q a_0))   \\
  & = \,  
  \sum_{j=0}^{n-2} \omega_q^j a_{n-j} + \omega_q^{n-1} (a_1 + \omega_q a_0) \\ 
 & = \, \sum_{j=0}^{n-2} \omega_q^j a_{n-j} + \omega_q^{n-1} a_1 + \omega_q^n a_0 
 \, = \, 
 \sum_{j=0}^{n} \omega_q^j a_{n-j} = S_q(a). 
\end{align*} 
\qed

\subsection*{(Left) Stacking the (Left) Trimmings for Binomial Summing Tests}

It would be remiss to mention the most well known summing tests, those 
that follow from the binomial identity. We wish to derive 
the binomial tests from an adapted form of Zbikowski's tests that trim from the 
left instead of the right, further solidifying the unification of trimming 
and summing tests. 

The binomial tests are developed by applying the binomial theorem to the 
standard expression for $a$ modulo $q$,  
$$ 
a \, \, = \, \, \sum_{j=0}^n {(q+(10-q))^j} a_j \, \, 
 \equiv \, \, \sum_{j=0}^n (10-q)^j a_j =: B_{q}(a) \mod q. 
$$
The well known tests for $9$ and $11$ are $B_9$ and $B_{11}$ 
and are usually motivated in this fashion. The binomial test for $7$ is  
$B_7(a) = {\sum_{j=0}^n }{3^j} {a_j}$ and for, say,  $39$ it is 
$B_{39}(a) = {\sum_{j=0}^n }{(-29)^j} {a_j}$. We claim that these tests 
can be developed via a recursive trimming and stacking procedure akin to 
the derivation of Khare's tests from Zbikowski's.  

On our main example, testing if $7$ divides $31284$, notice that 
we can rewrite $32184$ as $10^3( (7+3) \cdot 3 + 2) + 184$. The term in brackets is regarded 
as a non-traditional coefficient of $10^3$ just as we did in stacking (on the right) earlier. 
For testing divisibility by $7$ we can cast off the $7$ in the bracketed term before distributing, 
so $7$ divides $32184$ if and only if $7$ divides $10^3( (3) \cdot 3 + 2) + 184 = (11)184$. 
As before, this last number might be how we would write a check, writing the integer longhand 
as ``eleven thousand, one hundred, and eighty four.'' 

Repeating again, 
$(11)184 = 10^2(11(7+3) + 1) + 84$ reduces to $10^2(11(3) + 1) + 84 = (34)84$, 
or ``thirty four hundred and eighty four''. Repeating once more, 
$(34)84$ reduces to $(3 \cdot 34 + 8) 4 = (110)4$ which, repeating again, reduces to $(334)$. 
In other words, $7$ divides $32184$ if it divides $334$. We can repeat this process again on 
$334$ itself, should we wish, and it would equal $3^2 \cdot 3 + 3^1 \cdot 3 + 4 = 40$. We can 
conclude that $7$ does not divide $32184$. 

The example motivates an adapted version of Zbikowski's tests $T_q$ and the Stack function, 
which we will call \textit{left trim}, $LT_q$ and \textit{left stack}, LStack. 
It is immediate that 
$$
LT_q (a) := 10^{n-1}(10-q)a_n + a_[n-1,0]
$$ 
is a test for $q$ and that 
$$
\textup{LStack}(10^{n-1}(10-q)a_n + 10^{n-1} a_{n-1}+ a_{[n-2,0]}) = ( (10-q) a_n + a_{n-1})a_{[n-2,0]}
$$ 
provides the same flexibility that the original Stack function provided. Here is a more careful 
presentation of our main example with this notation. 

\vspace{.1in}

${\red 3}2184 \, { \xlongrightarrow{\, \, \, \, LT_7 \, \, \, \, }} \, \, 10^3 \cdot {\red (3 \cdot 3)} + 2184 \, $

\hspace{.35in}
$
{\red {\xlongequal{\bf LStack}} } \,   {\red (3 \cdot 3 + 2)}184$

\hspace{.35in}
$ \, { \xlongrightarrow{\, \, \, \, LT_7 \, \, \, \, }} \, \, 10^2 \cdot {\red 3 \cdot (3 \cdot 3 + 2) } + 184 \, $

\hspace{.35in}
$
{\red {\xlongequal{\bf LStack}} } \,   {\red (3^2 \cdot 3 + 3 \cdot 2 + 1)}84$

\hspace{.35in}
$ \, { \xlongrightarrow{\, \, \, \, LT_7 \, \, \, \, }} \, \, 10^1 \cdot {\red 3 \cdot (3^2 \cdot 3 + 3 \cdot 2 + 1) } + 84 \, $

\hspace{.35in}
$
{\red {\xlongequal{\bf LStack}} } \,  {\red (3^3 \cdot 3 + 3^2 \cdot 2 + 3 \cdot 1 + 8) } + 4$ 

\hspace{.35in}
$\, { \xlongrightarrow{\, \, \, \, LT_7 \, \, \, \, }} \, \, 10^0 \cdot {\red 3 \cdot (3^3 \cdot 3 + 3^2 \cdot 2 + 3 \cdot 1 + 8) } + 4 \, $

\hspace{.35in}
$
{\red {\xlongequal{\bf LStack}} } \,   {\red (3^4 \cdot 3 + 3^3 \cdot 2 + 3^2 \cdot 1 + 3 \cdot 8 + 4) } = {\red 334}$.

\vspace{.1in}

\begin{theorem} \label{the:binsum} 
The binomial test $B_{q}(a)$ equals $(\textup{LStack} \circ LT_q)^{n} (a)$.  
\end{theorem}

The proof is very similar to that of Theorem~\ref{the:sum}, inducting on the 
length of $a$ and trimming and stacking on the left as we did previously on 
the right. We leave the details as an exercise. 

\subsection*{Closing Remarks} 

Starting with the test for $7$ and using only elementary tools, we reduced Zbikowski's tests to 
a single trimming test for all integers. From Zbikowski's tests we derived Khare's 
summing tests as well as the binomial tests, adding only a dash of the induction axiom to our basic 
divisibility criteria. The two families of divisibility tests, trimming and summing, 
are much closer than initially meets the eye. 

Khare's tests are vastly preferable to the binomial tests and, in practice, the trimming tests 
are superior to both summing tests. The weights in Khare's tests {\em scale down} the original divisor 
$q$ by a factor of $10$ or $10/3$ whereas the binomial tests have weights that 
are the {\em difference} of $q$ with $10$. For example, Khare's $S_{39}(a) = {\sum_{j=0}^n }{4^j} {a_j}$
is preferable to the binomial $B_{39}(a) = {\sum_{j=0}^n }{(-29)^j} {a_j}$. The practice of Zbikowski's 
trimming is better than both as it avoids the mental computation of high powers of $\omega_q$, 
relying only on multiplying the last digit of an integer $a$ by $\omega_q$ followed by a straightforward 
subtraction and then recursively repeating this procedure. 

For Zbikowski's tests $a \, \equiv \, 10^{-1} T_q(a) \mod q$ and 
since stacking changes the representation of the number $a$ but not $a$ itself, then 
$a \equiv \omega_q T_q(a)$ and $a \equiv \omega_q^n S_q(a) \mod q$ whenever $a$ has length $n+1$. 
In contrast, part of the appeal of the binomial tests $B_q$ is its preservation of remainders. 

Khare also generalized the base $b=10$ to tests $S_q$ for $q$ in any base $b$. 
If $q$ and $b$ are co-prime then the $\omega_q$ term is precisely the least residue of 
$\omega_q \equiv b^{-1} \mod q$ and there are last-digits tests for all factors of $b$.  
Indeed, this article could be written for a general base $b$ and the results would 
hold as one would expect. 

Finally, while most tests are of the trimming and summing variety, there are tests 
that are not equivalent to those outlined here, like the Talmud test $\textup{Tal}_7$.  
Dickson \cite[Chapter XII]{Dickson} has many gems not discussed here and independently 
deriving each of them and understanding the some of the original sources would make for 
an excellent senior project. 

\subsection*{Acknowledgements}
I am grateful for conversations with 
Elizabeth Brown, Ezra ``Bud'' Brown, Brant Jones, Rachel Quinlan, and  
Jason Rosenhouse throughout my thinking and writing about divisibility tests. 
My student, Cameron Stopak suggested trimming from the left as an adaption to 
Zbikowski's trimming from the right. 
I was stuck for some time on how to jump from trimming to summing by elementary 
means and it was from playing shop with my children that I realized that 
``stacking pennies'' was exactly what was needed so thanks to e- and f- too.

\end{document}